\documentclass[12pt,reqno]{amsart}

\usepackage{amsmath,amsfonts,amsthm,amssymb,amsxtra}
\usepackage{bbm} % to get \1
\usepackage{hyperref}	% to get hyperlinks in equations and references
\newtheorem{theorem}{Theorem}%[section]

\newtheorem{lemma}[theorem]{Lemma}
\newtheorem{corollary}[theorem]{Corollary}

\theoremstyle{definition}

\newtheorem{definition}[theorem]{Definition}

\theoremstyle{remark}

\newtheorem{remark}[theorem]{Remark}

\numberwithin{equation}{section}
\numberwithin{theorem}{section}

\newcommand{\C}{\mathbb{C}}

\renewcommand{\epsilon}{\varepsilon}

\newcommand{\N}{\mathbb{N}}

\renewcommand{\phi}{\varphi}
\newcommand{\R}{\mathbb{R}}

\newcommand{\Z}{\mathbb{Z}}

\begin{document}

\title[Symmetrization inequalities on one-dimensional integer lattice ]{Symmetrization inequalities on one-dimensional integer lattice}

\author{By Shubham Gupta\\}
\address[Shubham Gupta]{Department of Mathematics, Imperial College London, 180 Queen’s Gate, London,
SW7 2AZ, United Kingdom.}
\email{s.gupta19@imperial.ac.uk}

\begin{abstract}
In this paper, we develop a theory of symmetrization on the one dimensional integer lattice. More precisely, we associate a radially decreasing function $u^*$ with a function $u$ defined on the integers and prove the corresponding Polya-Szeg\"{o} inequality. Along the way we also prove the weighted Polya-Szeg\"{o} inequality for the decreasing rearrangement on the $\emph{half-line}$, i.e., non-negative integers. As a consequence, we prove the discrete weighted Hardy's inequality with the weight $n^\alpha$ for $1 < \alpha \leq 2$. 
\end{abstract}

\maketitle

\section{Introduction}
\let\thefootnote\relax\footnotetext{Keywords: Decreasing rearrangement, Fourier rearrangement, symmetric-decreasing rearrangement, Polya-Szeg\"{o} inequality.\vspace{3pt}

2020 Mathematics Subject Classification: Primary: 39B62; Secondary: 42A16.
}
In this paper, the main inequality we are interested in is the discrete analogue of the \emph{Polya-Szeg\"{o} inequality} which states that
\begin{equation}\label{1.1}
    \int_{\R^d} |\nabla u|^2 dx \geq \int_{\R^d} |\nabla u^*|^2 dx,
\end{equation}
where $u^*$ is the \emph{symmetric-decreasing rearrangement} of the function $u$, defined by
\begin{equation}\label{1.2}
    u^*(x) := \int_0^\infty \chi_{\{|u|>t\}^*}(x) dt,
\end{equation}
where, given a set $\Omega$, $\Omega^*$ is the ball centered at origin whose measure is same as the measure of $\Omega \subseteq \R^d$. See \cite [chapter 3]{lieb2001analysis} for an introductory text on the subject of rearrangements and \cite[chapter 2]{kesavan2006symmetrization} for proofs of inequality \eqref{1.1}. It is important to note that $u^*$ is a radially decreasing function which is \emph{equimeasurable} with the function $|u|$, i.e., level sets of $|u|$ and $u^*$ have the same measure. Due to its vast applications, inequality \eqref{1.1} has been extended in various directions since the mid 20$^{th}$ century, we only cite a few papers here: \cite{ziemer1988minimal},\cite{brock1999weighted}, \cite{trombetti1997convex},\cite{cianchi2002functions}.\\

In Euclidean spaces $\R^d$, symmetrization inequalities including \eqref{1.1} are a standard and powerful tool to establish the symmetry of optimizers of many variational problems \cite{lenzmann2021sharp}, \cite{lieb1977existence}, \cite{henrot2006extremum}. In the past, they have also been successfully used to provide the sharp constant in various important functional inequalities(Sobolev inequality, Cafferreli-Kohn-Nirenberg inequality, etc.) \cite{talenti1976best}, \cite{alvino2017some}, \cite{alvino2019weighted}, \cite{lieb2002sharp}, \cite{moser1971sharp}. However, this very important analytic tool seems to be missing in the discrete setting, for instance, integer lattices: $\Z^d$, and graphs in general. Till date, very little is known about the optimizers of discrete versions of many important functional inequalities. For instance, existence of optimizers of Sobolev inequality on lattices was established very recently in a paper by Bobo Hua and Ruowei L in 2021 \cite{hua2021existence}. Another important example is that of Hardy's inequality. Nothing is known about the sharp constant and optimizers of  Hardy's inequality on lattices, in fact, to our best knowledge, \cite{kapitanski2016continuous}is the only paper where an explicit constant was computed for Hardy's inequality on $d$-dimensional lattices. We hope that by developing a theory of discrete symmetrization, one might be able to answers these questions. \\

In this work, we develop a theory of symmetrization on the one dimensional lattice $\Z$. In particular we are interested in proving the following discrete analogue of \eqref{1.1}
\begin{equation}\label{1.3}
    \sum_{n \in \Z}|u(n)-u(n-1)|^2 \geq \sum_{n \in \Z} |u^*(n)-u^*(n-1)|^2,
\end{equation}
for some suitable rearrangement $u^*$ of the function $u$. We hope that this work might shed some light on how to develop a theory on general $d-$dimensional lattices. Before getting into the details of our results, we would like to mention previous works that have been done in this direction. The oldest results we could find dates back to the early 20th century. These results are compiled in chapter X of the book by Hardy, Littlewood and Polya \cite{hardy1952inequalities}. The chapter contains important rearrangement inequalities like the \emph{Hardy-Littlewood inequality} and the  \emph{Riesz-rearrangement inequality} for compactly supported functions defined on integers. In late 90's Alexander R. Pruss wrote a couple of papers on Riesz-rearrangement inequality on general graphs and their applications\cite{pruss1998discrete},\cite{Pruss_1996}(amongst many important things, they proved \emph{Faber-Krahn} inequality on regular trees). In 2010 \cite{hajaiej2010rearrangement}, Hichem Hajaiej proved Polya-Szeg\"{o} inequality \eqref{1.3} on integers using the idea of polarization. The idea of polarization has been very effective in proving various symmetrization inequalities on euclidean spaces, see for instance \cite{burchard2006rearrangement}, \cite{brock2000approach}, \cite{brock2000general}, \cite{draghici2005rearrangement}. It seems, no further progress has been made in the direction of discrete symmetrization after Hajaiej's paper. \\

The major problem in discretizing the definiton \eqref{1.2} on $\Z^d$ is that given a set $\Omega$ in $\Z^d$ it is not always possible to find a ball in the lattice whose measure coincides with that of $\Omega$. Due to this reason, one is bound to lose the radiality property of $u^*$, which is important for applications. In all the works mentioned above, discrete rearragements only has the monotonicity property and lacks the radiality property. For instance, in \cite{hajaiej2010rearrangement}, for a function $u: \Z \rightarrow \R$, the author defines $u^*$ as follows: $u^*(0)$ is the largest value of $|u|$, $u^*(1)$ is the second largest value of $|u|$, $u^*(-1)$ is the third largest value of $|u|$, and so on. Essentially $u^*$ is a decreasing function with respect to some \emph{spiral-like ordering} of the set of integers. One of the novel features of this paper is that we define a rearrangement $u^*$ on integers, which is both radial and monotonic, and we prove the corresponding Polya-Szeg\"{o} principle \eqref{1.3}. This is done in three steps: in section \ref{sec:decreasing} we define the standard decreasing rearrangement on the \emph{half-line}(i.e., on non-negative integers) and prove the weighted version of the Polya-Szeg\"{o} inequality \eqref{1.3}. We would like to mention that although the standard Polya-Szeg\"{o} inequality for decreasing rearrangement is already known(see \cite{hajaiej2010rearrangement}, for instance), the weighted version seems to be missing from current literature. Moreover, our proof of the unweighted case via the discrete co-area formula is also new. Secondly, we define \emph{Fourier rearrangement} in section \ref{sec:fourier}, which associates to every function $u \in l^2(\Z)$, a radial function $u^\#$. We prove Polya-Szeg\"{o} principle for the Fourier rearrangement, furthermore, we prove that Polya-Szeg\"{o} holds true even when one replaces first order derviatives with higher order derivatives. A similar type of Fourier rearrangement procedure in $\R^d$ was carried out recently by Lenzmann and Sok \cite{lenzmann2021sharp}. In that paper, they proved Polya-Szeg\"{o} for higher order operators and used it to prove radial symmetry of various higher order variational problems. Our Fourier rearrangement method can be thought of as a discrete analogue of their work. Finally, in section \ref{sec:symmetricdecreasing}, we combine both decreasing and Fourier rearrangements to define the \emph{symmetric-decreasing rearrangement} $u^*$, of the function $u$, which is always a radially decreasing function and we prove the corresponding Polya-Szeg\"{o} inequality for this rearrangement.\\

In the last section \ref{sec:applications}, we apply the weighted Polya-Szeg\"{o} inequality proved in section \ref{sec:decreasing} to prove discrete Hardy's inequality on non-negative integers with weights $n^\alpha$ for $1 < \alpha \leq 2$, refer \cite{gupta2021discrete} and \cite{gupta2021one} for the proof of this inequality for $ \alpha \in [0,1) \cup [5, \infty)$ and $\alpha \in 2\N$ respectively.

\section{Decreasing rearrangement on the half line}\label{sec:decreasing}
In this section, we will be dealing with functions which vanish at infinity in the following sense: Let $\Z^+$ be the set of all non-negative integers and let $u: \mathbb{Z}^+ \rightarrow \C$ be a function. We say that function $u$ $\emph{vanishes at infinity}$ if all level sets of $|u|$ are finite, i.e., $\{x: |u(x)| > t\}$ is finite for all $t > 0$. \\

Let $A$ be a finite subset of $\mathbb{Z}^+$ such that $|A| = k$. Then we define its \emph{symmetric rearrangement} $A^*$ as the set of all non-negative integers which are less than or equal to $k-1$. Thus, $A^*$ is a ball in $\mathbb{Z}^+$ centered at origin whose size is same as the size of $A$. \\

Let $u: \mathbb{Z}^+ \rightarrow \C$ be a function vanishing at infinity. We define its \emph{decreasing rearrangement} $\widetilde{u}$ as:
\begin{equation}
    \widetilde{u}(n) := \int_0^\infty \chi_{\{|u|>t\}^*}(n) dt.
\end{equation}
Some useful properties of $\widetilde{u}$:
\begin{enumerate}
    \item $\widetilde{u}$ is always non-negative. \\
    \item The level sets of $\widetilde{u}$ are rearrangements of level sets of $|u|$, i.e., 
    \begin{equation}
        \{n: \widetilde{u} > t\} = \{n: |u|>t\}^*,
    \end{equation}
    for all $t >0$. An important consequence of this is that level sets of $\widetilde{u}$ and $|u|$ have the same size. \\
    \item It is easy to conclude that $\widetilde{u}(0)$ is the largest value of $|u|$, $\widetilde{u}(1)$ is the second largest value of $|u|$, $\widetilde{u}(2)$ is the third largest value of $|u|$ and so on. In particular, $\widetilde{u}$ is a decreasing function of $n$.\\
    \item Equimeasurability of levels sets of $\widetilde{u}$ and $|u|$ along with layer cake representation gives
    \begin{equation}\label{2.3}
        \sum_{n \in \mathbb{Z}^+}|u|^p = \sum_{n \in \mathbb{Z}^+}|\widetilde{u}|^p,
    \end{equation}
    for $p \geq 1$.\\
    \item Let $\Phi: \mathbb{R^+} \rightarrow\mathbb{R^+}$ be a bijective map. Then 
    \begin{equation}
        \widetilde{\Phi(|u|)} = \Phi(\widetilde{u}).
    \end{equation}
    \item The rearrangement is \emph{order preserving}, i.e, if $u$ and $v$ are two non-negative functions such that $u(n) \leq v(n)$ then $\widetilde{u}(n) \leq \widetilde{v}(n)$. This follows from the fact that $u(n) \leq v(n)$ is equivalent to $\{u>t\} \subseteq \{v>t\}$.\\
    \item (\emph{Hardy-Littlewood inequality}) Let $u$ and $v$ be non-negative functions vanishing at infinity. Then we have 
    \begin{equation}\label{2.5}
        \sum_{n \in \mathbb{Z}^+}u(n)v(n) \leq \sum_{n \in \mathbb{Z}^+}\widetilde{u}(n)\widetilde{v}(n). 
    \end{equation}
    Using layer-cake representation of functions $u, v$ and monotone convergence theorem we get
    \begin{align*}
        \sum_{n \in \mathbb{Z}^+}u(n)v(n) = \int_0^\infty \int_0^\infty \sum_{n \in \Z^+} \chi_{\{u>t\}}(n) \chi_{\{v>s\}}(n) dt ds. 
    \end{align*}
    Therefore it is enough to prove \eqref{2.5} when $u$ and $v$ are characteristic functions, i.e., $u = \chi_A$ and $v = \chi_B$ for finite subsets $A$ and $B$ of $\mathbb{Z}^+$. For this choice of $u$ and $v$ inequality \eqref{2.5} is equivalent to proving $|A \cap B| \leq |A^* \cap B^*|$. Without loss of generality we can assume that $|A|\leq |B|$, which implies $A^* \subseteq B^*$. Then $|A^* \cap B^*| = |A^*|= |A|\geq |A \cap B|$. This completes the proof of the Hardy-Littlewood inequality \eqref{2.5}.\\
    
    \item $l^p(\Z^+)$ distance decreases under decreasing rearrangement, i.e., 
    \begin{equation}\label{2.6}
        \sum_{n \in \Z^+}|\widetilde{u}(n) -\widetilde{v}(n)|^p \leq \sum_{n \in \Z^+} |u(n)-v(n)|^p,
    \end{equation}
    for all $u, v$ non-negative functions vanishing at infinity and $p \geq 1$. For $p=2$, inequality \eqref{2.6} follows directly from \eqref{2.3} and Hardy-Littlewood inequality \eqref{2.5}. For the proof of the general case $p\geq 1$, we refer \cite[theorem 3.5, chapter 3]{lieb2001analysis}, where inequality \eqref{2.6} was proved for rearrangements on the real line.  \\
\end{enumerate}

Next we will prove the main theorem of this section, a \emph{weighted Polya-Szeg\"{o} inequality}. We would like to mention that the Polya-Szeg\"{o} inequality on the integers has been considered in the past \cite{hardy1952inequalities, pruss1998discrete, hajaiej2010rearrangement}, but to the author's best knowledge, the Polya-Szeg\"{o} principle with weights has not been considered before. The major tool in our proof will be the discrete analogue of the co-area formula, which is the content of our next lemma. We will state the co-area formula in the general setting of graphs.
\begin{lemma}[Discrete co-area formula]
Let $G =(V,E)$ be a graph, where $V$ is a countable set of vertices, and $E$ is a symmetric relation on set $V$. We also assume that $G$ is locally finite, i.e., each vertex of $G$ has finite degree. Let $u: V(G) \rightarrow \mathbb{R}$ be a non-negative function. For $p \geq 1$ we have,
\begin{equation}\label{2.7}
    \sum_{x \sim y}|u(x)-u(y)|^p = \int_0^\infty \sum_{(x,y) \in E(\{u>t\}, \{u>t\}^c)}|u(x)-u(y)|^{p-1} dt,
\end{equation}
where $E(X,Y)$ denotes the set of edges between subsets $X$ and $Y$ of graph $G$, and $x \sim y$ means vertices $x$ and $y$ are connected by an edge.
\end{lemma}
\begin{proof}
Let $A: V \times V \rightarrow \R$ be a map defined as:
\begin{align*}
    A(x,y) := 
    \begin{cases}
    &1 \hspace{19pt} \text{if} \hspace{5pt} x \sim y\\
    &0 \hspace{19pt} \text{if} \hspace{5pt} x \nsim y
    \end{cases}
\end{align*}
Consider, 
\begin{align*}
    \sum_{x \sim y}|u(x)-u(y)|^p &= \sum_{x \in V} \sum_{y \in V} A(x,y)|u(x)-u(y)|^p\\
    &= \sum_{u(x) \leq u(y)}A(x,y)|u(x)-u(y)|^{p-1}(u(y)-u(x)) \\
    &+ \sum_{u(x) > u(y)}A(x,y)|u(x)-u(y)|^{p-1}(u(x)-u(y))\\
    &= \int_0^\infty  \sum_{u(x) \leq u(y)}A(x,y)|u(x)-u(y)|^{p-1} \chi_{[u(x),u(y))}(t) dt\\
    &+ \int_0^\infty \sum_{u(x) > u(y)}A(x,y)|u(x)-u(y)|^{p-1}  \chi_{[u(y),u(x))}(t) dt.
\end{align*}
Finally, using the symmetry of the integrand with respect to the variables $x$ and $y$ we obtain,
\begin{align*}
    \sum_{x \sim y}|u(x)-u(y)|^p = \int_0^\infty \sum_{(x,y) \in E(\{u>t\}, \{u> t\}^c)}|u(x)-u(y)|^{p-1} dt. 
\end{align*}
\end{proof}

\begin{theorem}[Weighted Polya-Szeg\"{o} inequality]\label{thm2}
Let $w: \Z^+ \rightarrow \R$ be a non-negative increasing function. Let $u: \Z^+ \rightarrow \C$ be a function which is vanishing at infinity. We have, 
\begin{equation}\label{2.8}
    \sum_{n \in \Z^+} |u(n)-u(n+1)|^p w(n) \geq \sum_{n \in \Z^+} |\widetilde{u}(n)-\widetilde{u}(n+1)|^p w(n), 
\end{equation}
for  all $p \geq 1$. Moreover, if $w(n) >0$, if $u$ produces equality in \eqref{2.8}, then $|u| = \widetilde{u}$. In particular, $|u|$ is a decreasing function.
\end{theorem}

\begin{remark}
Choosing $w(n)=|n|^\alpha$ for non-negative $\alpha$ in theorem \ref{thm2} gives 
\begin{equation}\label{2.9}
    \sum_{n \in \Z^+}|u(n)-u(n+1)|^p |n|^\alpha \geq \sum_{n \in \Z^+}|\widetilde{u}(n)-\widetilde{u}(n+1)|^p |n|^\alpha,
\end{equation}
which is the discrete analogue of Corollary 8.1 of the paper \cite{alvino2017some}.
\end{remark}

\begin{proof}
The proof goes via the discrete co-area formula. We interpret $\Z^+$ as a graph with $V = \Z^+$. Two points $x$ and $y$ are connected by an edge iff $|x-y| = 1$. Using the co-area formula \eqref{2.7}, we get, 
\begin{equation}\label{2.10}
    \sum_{x \sim y}|u(x)-u(y)|^p w((x+y-1)/2) = \int_0^\infty \sum_{(x, y) \in \partial(\{u>t\})} |u(x)-u(y)|^{p-1}w((x+y-1)/2) dt,
\end{equation}
where $\partial(\{u>t\}) := E(\{u>t\},\{u>t\}^c)$.\\

Using the reverse triangle inequality, we get $\sum \limits_{n \in \Z^+}|u(n)-u(n+1)|^p w(n) \geq \sum \limits_{n \in \Z^+}||u|(n)-|u|(n+1)|^p w(n)$. Therefore, without loss of generality, we can assume that $u$ is non-negative. Since $u$ vanishes at infinity, we can arrange the values $u$ takes in decreasing order: $t_1 > t_2 > t_3 ....$ and so on. Fix $ t_{i+1} \leq t < t_i$. Let $k$ be the size of the level set of $u$, i.e., $|\{u>t\}| = k$. Then it is easy to see that there will exist $x_1 \geq k-1$ such that $x_1 \in \{u>t\}$ and $x_1+1 \in \{u>t\}^c$. This is a consequence of $u$ vanishing at infinity. This implies that,
\begin{equation}\label{2.11}
    \begin{split}
        \frac{1}{2}\sum_{(x, y) \in \partial(\{u>t\})} |u(x)-u(y)|^{p-1}w((x+y-1)/2) & \geq |u(x_1)-u(x_1+1)|^{p-1}w(x_1) \\
        & \geq |t_i -t_{i+1}|^{p-1}w(k-1).   
    \end{split}
\end{equation}
Further, notice that $\{\widetilde{u}>t\} = [0, k-1]\cap \Z^+$ and
\begin{equation}\label{2.12}
    \frac{1}{2}\sum_{(x, y) \in \partial(\{\widetilde{u}>t\})} |\widetilde{u}(x)-\widetilde{u}(y)|^{p-1}w((x+y-1)/2)
    = |t_i - t_{i+1}|^{p-1}w(k-1). 
\end{equation}
Equation \eqref{2.11} and \eqref{2.12} imply 
\begin{equation}\label{2.13}
    \sum_{(x, y) \in \partial(\{u>t\})} |u(x)-u(y)|^{p-1}w((x+y-1)/2) \geq \sum_{(x, y) \in \partial(\{\widetilde{u}>t\})} |\widetilde{u}(x)-\widetilde{u}(y)|^{p-1}w((x+y-1)/2), 
\end{equation}
for all $t>0$. Equation \eqref{2.13} along with the coarea formula \eqref{2.10} completes the proof.\\

Next we will characterize those functions which attain equality in \eqref{2.8}. Let $u$ be a non-negative function vanishing at infinity which produces equality in \eqref{2.8}. We have
\begin{equation}\label{2.14}
    \int_0^\infty \sum_{(x, y) \in \partial(\{u>t\})} |u(x)-u(y)|^{p-1}w((x+y-1)/2) dt = \int_0^\infty \sum_{(x, y) \in \partial(\{\widetilde{u}>t\})} |\widetilde{u}(x)-\widetilde{u}(y)|^{p-1}w((x+y-1)/2) dt.
\end{equation}
Identity \eqref{2.14} along with \eqref{2.13} implies that
\begin{equation}\label{2.15}
    \sum_{(x, y) \in \partial(\{u>t\})} |u(x)-u(y)|^{p-1}w((x+y-1)/2) = \sum_{(x, y) \in \partial(\{\widetilde{u}>t\})} |\widetilde{u}(x)-\widetilde{u}(y)|^{p-1}w((x+y-1)/2),
\end{equation}
for a.e. $t > 0$. Let $t_{i+1} \leq  t < t_i $, where $t_1>t_2>t_3>...$ are the values of function $u$ arranged in decreasing order. Let $|\{u>t\}| = k$. Then we claim that there cannot exist $x \geq k$ such that $u(x) > t$. We will prove this via contradiction. Let us assume that such an $x$ exists. It is easy to see that there will exist $x_1 \geq k$ and $x_2 \leq x$ such that $x_1, x_2 \in \{u>t\}$ and $x_1+ 1, x_2 -1 \in \{u>t\}^c$. This is a consequence of the fact that the level set of $u$ has size $k$, which is a finite number. This immediately implies that
\begin{align*}
   \frac{1}{2}\sum_{(x, y) \in \partial(\{u>t\})} |u(x)-u(y)|^{p-1}w((x+y-1)/2) \geq &|u(x_1)-u(x_1 +1)|^{p-1}w(x_1) \\
   & + |u(x_2)-u(x_2-1)|^{p-1}w(x_2-1)\\
   & > |t_i - t_{i+1}|^{p-1}w(k-1). 
\end{align*}
The last step uses the fact that $w$ is a positive function. This contradicts identity \eqref{2.15}. Therefore $\{u>t\} = \{\widetilde{u}>t\}$ for a.e. $t>0$. Finally, using layer-cake representation for $u$ and $\widetilde{u}$, we get $u(x) = \widetilde{u}(x)$. \\

Let us assume that function $u$ produces equality in \eqref{2.8}, then $|u|$ will also produce equality in \eqref{2.8} which would imply that $|u| = \widetilde{u}$.  
\end{proof}

\begin{remark}
We would like to point out that \eqref{2.8} is not true in general for higher order operators. Consider the following counter example: let $u(0) = u(2) = \alpha$ and $u(1) = \alpha + \delta$ and define $u = 0$ rest everywhere. Then 
\begin{align*}
    \sum_{n \in \Z^+}|\Delta u|^2  = \alpha^2 + 5\delta^2 + (\alpha-\delta)^2 \hspace{5pt} \text{and} \hspace{5pt} \sum_{n \in \Z^+}|\Delta \widetilde{u}|^2 = 2(\alpha^2 + \delta^2)   
\end{align*}
where we define $\Delta u(n) := 2u(n)-u(n-1)- u(n+1)$ for $n \geq 1$ and $\Delta u(0):=  u(0)-u(1)$. It is easy to check that $\sum \limits_{n \in \Z^+}|\Delta \widetilde{u}|^2 > \sum \limits_{n \in \Z^+} |\Delta u|^2$
for $0 < \delta \leq \alpha/2$.
\end{remark}

\section{Fourier rearrangement}\label{sec:fourier}
In this section, we will introduce another rearrangement called the  \emph{Fourier rearrangement} on $\Z$ and prove the corresponding Polya-Szeg\"{o} principle. As we will see, under this rearrangement, norms of higher order operators decrease as well, something which was missing in the decreasing rearrangement introduced in the last section. But on the downside, the description of fourier rearrangement of a function is not as straightforward as that of the decreasing rearrangement. Many basic yet important properties like equimeasurability of level sets fail to hold in this setting. \\

Before mentioning the details, we would like to mention \cite{lenzmann2021sharp}, where an analogous theory of Fourier rearrangement on $\R^d$ has been developed and applied. \\

Let $u \in l^2(\mathbb{Z})$. Then we define its Fourier transform $F(u) \in L^2((-\pi, \pi))$ by 
\begin{align*}
    F(u)(x) := (2\pi)^{-\frac{1}{2}}\sum_{n \in \mathbb{Z}} u(n) e^{-inx} \hspace{19pt} x \in(-\pi, \pi),
\end{align*}

and the inverse of $F$ is given by
\begin{equation}
    F^{-1}(u)(n) = (2\pi)^{-\frac{1}{2}} \int_{-\pi}^\pi u(x) e^{inx} dx
\end{equation}
We define the \emph{Fourier rearrangement} $u^\#$ \emph{of the function} $u \in l^2(\mathbb{Z})$ as 
\begin{align*}
    u^\# := F^{-1}((F(u))^*),
\end{align*}
where $f^*$ denotes the symmetric decreasing rearrangement of the function $f$ on the real line (see chapter 3 of \cite{lieb2001analysis}).\\\\
Some basic properties of $u^\#$:
\begin{enumerate}
    \item $u^\#$ is radial and real valued. To prove this fact consider,
    \begin{align*}
        u^\#(-n) &= (2\pi)^{\frac{1}{2}} \int_{-\pi}^\pi F(u)^*(x) e^{-inx}dx = (2\pi)^{\frac{1}{2}}  \int_{-\pi}^\pi F(u)^*(-x) e^{inx}dx = u^\#(n).\\
        \overline{u^\#(n)} &= (2\pi)^{-\frac{1}{2}}\int_{-\pi}^\pi \overline{F(u)^*(x)} e^{-inx}dx = (2\pi)^{\frac{1}{2}} \int_{-\pi}^\pi F(u)^*(-x) e^{inx}dx = u^\#(n).
    \end{align*}
    \item $|u^\#(n)| \leq |u^\#(0)|$, i.e., $|u^\#|$ takes its maximum value at origin. This follows directly from the formula for the inverse of fourier transform:
    \begin{align*}
        |u^\#(n)| \leq (2\pi)^{-\frac{1}{2}} \int_{-\pi}^\pi |F(u)^*| dx = u^\#(0). 
    \end{align*}
    \item Using Parseval's identity and equimeasurability of symmetric decreasing rearrangement we get,   
    \begin{equation}\label{3.2}
        \sum_{n \in \mathbb{Z}} |u|^2 = \sum_{n \in \mathbb{Z}}|u^\#|^2. 
    \end{equation}
    \item  We have the following \emph{Hardy-Littlewood inequality}:
    \begin{equation}\label{3.3}
        \Big|\sum_{n \in \mathbb{Z}} u(n) \overline{v}(n)\Big| \leq \sum_{n \in \mathbb{Z}} u^\#(n) v^\#(n). 
    \end{equation}
    Inequality \eqref{3.3} follows from Parseval's identity and the Hardy-Littlewood inequality for symmetric decreasing rearrangement.
    \item $l^2(\Z)$ distance decreases under Fourier rearrangement, i.e.,
    \begin{align*}
        \sum_{n \in \mathbb{Z}}|u^\#-v^\#|^2 \leq \sum_{n \in \mathbb{Z}} |u - v|^2.
    \end{align*}
    for all $u, v \in l^2(\Z)$. It again follows from Parseval's identity and the $L^2(\R)$ contraction property of symmetric decreasing rearrangement.\\
\end{enumerate}
Before stating the main result, we will define the discrete analogue of the first and second order derivatives on integers.

\begin{definition}
Let $u: \Z \rightarrow \C$ be a function defined on integers, then its \emph{first and second order derivative} are defined by
\begin{equation}
    Du(n) := u(n)-u(n-1),
\end{equation}
and
\begin{equation}
    \Delta u:= 2u(n)-u(n-1)-u(n+1)
\end{equation}
respectively.
\end{definition}

\begin{theorem}\label{thm3.2}

Let $u \in l^2(\mathbb{Z})$ and $k \geq 0$ then we have
\begin{equation}\label{3.6}
    \sum_{n \in \mathbb{Z}} |\Delta^k u|^2 \geq \sum_{n \in \mathbb{Z}} |\Delta^k u^\#|^2, 
\end{equation}
and 
\begin{equation}\label{3.7}
    \sum_{n \in \mathbb{Z}} |D\Delta^k u|^2 \geq \sum_{n \in \mathbb{Z}} |D\Delta^k u^\#|^2.
\end{equation}
Moreover equality holds in \eqref{3.6} or \eqref{3.7} if and only if $|F(u)(x)| = F(u)^*(x)$ for a.e. $x \in (-\pi, \pi)$. 
\end{theorem}

\begin{remark}
Putting $k=0$ in inequality \eqref{3.7} gives the Polya-Szeg\"{o} inequality \eqref{1.3} for the Fourier rearrangement $u^\#$.
\end{remark}

\begin{remark}
We would like to mention that theorem \ref{thm3.2} is true for general class of operators. Let $T: l^2(\Z) \rightarrow l^2(\Z)$ be an operator. Assume that that there exists a measurable function $\omega : \R \rightarrow \C$ such that $|\omega(x)|$ is radial, strictly increasing with respect to $|x|$ and 
\begin{align*}
    F(T(u))(x) = \omega(x)F(u)(x),
\end{align*}
for all $x \in (-\pi, \pi)$. Then we have
\begin{align*}
    \sum_{n \in \Z} |Tu|^2 \geq \sum_{n \in \Z}|Tu^\#|^2.
\end{align*}
It is easy to see that operators $\Delta^k$ and $D \Delta^k$ satisfy the above conditions. One can also put assumptions on the parameters $a,b,c$ for which the above conditions are satisfied by \emph{Jacobi operators}: $Ju(n) := au(n) + bu(n-1) + cu(n+1)$. 
\end{remark}

\begin{proof}
We begin by computing the Fourier transforms of $Du$ and $\Delta u$. 
\begin{align*}
    Du(n) = u(n)-u(n-1) &= (2\pi)^{-\frac{1}{2}}\int_{-\pi}^\pi F(u)(1-e^{-ix}) e^{inx} dx,\\
    \Delta u(n) = 2u(n) - u(n-1) - u(n+1) &= (2\pi)^{-\frac{1}{2}} \int_{-\pi}^\pi F(u)(2-e^{-ix}-e^{ix}) e^{inx} dx.
\end{align*}
Therefore we have 
\begin{align*}
    F(Du) = F(u)(1-e^{-ix}) \hspace{19pt} \text{and} \hspace{19pt} 
    F(\Delta u) = F(u)(2-e^{ix}-e^{-ix}).
\end{align*}
Using Parseval's identity and applying the above formulas iteratively we obtain
\begin{equation}\label{3.8}
    \sum_{n \in \mathbb{Z}} |\Delta^k u|^2 = \int_{-\pi}^\pi |F(\Delta^k u)|^2 dx = 4^{2k} \int_{-\pi}^\pi |F(u)|^2 \sin^{4k}(x/2) dx,
\end{equation}
and
\begin{equation}\label{3.9}
    \sum_{n \in \mathbb{Z}} |D\Delta^k u|^2 = \int_{-\pi}^\pi |F(D\Delta^k u)|^2 dx = 4^{2k+1} \int_{-\pi}^\pi |F(u)|^2 \sin^{4k+2}(x/2) dx.
\end{equation}
Therefore proving \eqref{3.6} and \eqref{3.7} reduces to showing
\begin{equation}\label{3.10}
    \int_{-\pi}^\pi |F(u)|^2 \sin^{4k}(x/2) dx \geq \int_{-\pi}^\pi |F(u)^*|^2 \sin^{4k}(x/2) dx,
\end{equation}
and 
\begin{equation}\label{3.11}
    \int_{-\pi}^\pi |F(u)|^2 \sin^{4k+2}(x/2) dx \geq \int_{-\pi}^\pi |F(u)^*|^2 \sin^{4k+2}(x/2) dx
\end{equation}
respectively. Next we will prove a general result which implies \eqref{3.10} and \eqref{3.11}.\\

\textbf{Claim:} Let $f$ be a non-negative measurable function vanishing at infinity and let $g$ be a non-negative measurable function which is  radially increasing. Then we have
\begin{equation}\label{3.12}
    \int_\mathbb{R} f(x)g(x)  dx \geq \int_\mathbb{R} f^*(x) g(x) dx.
\end{equation}
\begin{proof}
Using the Layer-cake representation formula proving \eqref{3.7} reduces to proving \eqref{3.7} for $f = \chi_\Omega$ for measurable set $\Omega $ of finite measure. Let $\Omega^*$ be the ball centered at origin such that $|\Omega^*| = |\Omega|$. Using the facts that $g$ is radially increasing and $|\Omega \setminus \Omega^*| = |\Omega^*\setminus \Omega|$ we obtain
\begin{align*}
    \int_{\Omega \setminus \Omega^*} g(x) dx \geq g(R)|\Omega \setminus \Omega^*| = g(R)|\Omega^* \setminus \Omega| \geq \int_{\Omega^*\setminus \Omega} g(x) dx,
\end{align*}
where $R$ is the radius of the ball $\Omega^*$. Now we have
\begin{align*}
    \int_\Omega g(x) dx = \int_{\Omega \setminus \Omega^*} g(x) dx + \int_{\Omega \cap \Omega^*} g(x) dx \geq \int_{\Omega^* \setminus \Omega} g(x) dx + \int_{\Omega^* \cap \Omega} g(x) dx = \int_{\Omega^*} g(x) dx.
\end{align*}
This prove inequality \eqref{3.12}. 
\end{proof}

Now coming back to inequality \eqref{3.10}, we define $f := |F(u)|^2 $ on the interval $(-\pi,\pi)$ and zero in the complement, and $g : = 4\sin^{4k}(x/2)$ in the interval $(-\pi, \pi)$ and extend $g(x)$ in a radial and strictly increasing way in the complement of $(-\pi, \pi)$. Now applying \eqref{3.12} for this choice of $f$ and $g$ yields \eqref{3.10}. One can similarly prove \eqref{3.11}, thereby completing the proof of inequalities \eqref{3.6} and \eqref{3.7}.\\

Next we will study the equality cases in \eqref{3.6} and \eqref{3.7}. For that, we will study the equality case in inequality \eqref{3.12}. Let $f$ be a non-negative function which produces equality in \eqref{3.12}. We further assume that $g(x)$ is a strictly increasing function, then we have for a.e. $t>0$,
\begin{equation}
    \int_{\{f>t\}} g(x) dx = \int_{\{f^*>t\}} g(x) dx.
\end{equation}
This implies that $\int_{\Omega \setminus \Omega^*} g(x) = g(R)|\Omega\setminus \Omega^*| $, with $\Omega := \{f>t\}$ and $R$ being the radius of the ball $\Omega^*$. Now we claim that $|\Omega \setminus \Omega^*| = 0$, if this is not true then the fact that $g$ is strictly increasing would imply that $\int_{\Omega \setminus \Omega^*} g(x) > g(R)|\Omega\setminus \Omega^*|$, which leads to a contradiction. Therefore, we have $|\Omega \setminus \Omega^*| =0$, which implies that 
\begin{align*}
    \chi_{\{f>t\}}(x) = \chi_{\{f^*>t\}}(x),
\end{align*}
for a.e. $(x,t) \in \R \times (0, \infty)$. Finally using the layer-cake representation for $f$ and $f^*$, we find that equality  in \eqref{3.12} implies that $f(x) = f^*(x)$ for a.e. $x \in \R$.\\

Let's get back to the equality case in \eqref{3.6}. Let $u$ be a function which produces equality in \eqref{3.6} then we have 
 
\begin{equation}
    \int_{-\pi}^\pi |F(u)|^2 \sin^{4k}(x/2) dx = \int_{-\pi}^\pi (|F(u)|^2)^* \sin^{4k}(x/2) dx,
\end{equation}
which implies that $|F(u)|^2 = (|F(u)|^2)^* = |F(u)^*|^2$. This gives us $|F(u)| = F(u)^*$ for a.e. $ x \in (-\pi, \pi)$. 

\end{proof}

\begin{remark}
We would like to point out that unlike the decreasing rearrangement, weighted Polya-Szeg\"{o} inequalities of the type \eqref{2.13} with power weights do not hold true in general for Fourier rearrangements. Let us give a example to support this statement. Let $u(0)=u(1)=\beta$ and $u$ vanishes everywhere else. Then $F(u) = 2\beta(2\pi)^{-\frac{1}{2}} e^{-ix/2}\cos(x/2)$ and its rearrangement is given by $F(u)^* = 2|\beta|(2\pi)^{-\frac{1}{2}}\cos(x/2)$. Finally, using the inversion formula we obtain
\begin{equation}
    u^\# = \frac{4|\beta|}{\pi}\frac{(-1)^n}{1-4n^2}.
\end{equation}
Consider,
\begin{align*}
    \sum_{n \in \mathbb{Z}}|u(n)-u(n-1)|^2|n|^\alpha = \beta^2\Big((1/2)^\alpha + (3/2)^\alpha \Big),
\end{align*}
and 
\begin{align*}
    \sum_{n \in \mathbb{Z}}|u^\#(n)-u^\#(n-1)|^2|n|^\alpha = \Big(\frac{4|\beta|}{\pi}\Big)^2\sum_{n \in \mathbb{Z}} |v(n)-v(n-1)|^2 |n|^\alpha,
\end{align*}
where $v(n):= \frac{(-1)^n}{1-4n^2}$.\\\\
It is easy to check that $|v(n)-v(n-1)|^2 |n|^\alpha$ grows as $|n|^{\alpha-4}$ as $|n| \rightarrow \infty$. So clearly, for $\alpha > 4$, we will have
\begin{align*}
    \sum_{n \in \mathbb{Z}}|u^\#(n)-u^\#(n-1)|^2|n|^\alpha >  \sum_{n \in \mathbb{Z}}|u(n)-u(n-1)|^2|n|^\alpha.   
\end{align*}
In fact we will have $\sum \limits_{n \in \mathbb{Z}}|u^\#(n)-u^\#(n-1)|^2|n-1/2|^\alpha = \infty$. Therefore we cannot have a real constant $c$ for which 
\begin{align*}
    c\sum_{n \in \mathbb{Z}}|u(n)-u(n-1)|^2|n|^\alpha \geq \sum_{n \in \mathbb{Z}}|u^\#(n)-u^\#(n-1)|^2|n|^\alpha
\end{align*}
holds.
\end{remark}

\begin{remark}
In this remark, we would illustrate how one can deduce formulae for some infinite sums using Fourier rearrangement. Let $u$ be the function as defined in the above remark, that is, $u(0)=u(1) = \beta$ and $u$ is zero everywhere else. Then we have
\begin{align*}
    u^\#(n) = \frac{4|\beta|}{\pi} \frac{(-1)^n}{1-4n^2},
\end{align*}
Then using the fact that $l^2(\Z)$ norm is preserved under fourier rearrangement \eqref{3.2} we get,
\begin{align*}
    2\beta^2 = \Big(\frac{4\beta}{\pi}\Big)^2 \sum_{n \in \Z} \frac{1}{(4n^2-1)^2},
\end{align*}
which implies,
\begin{equation}
    \sum_{n \in \Z} \frac{1}{(4n^2-1)^2} = \pi^2/8.
\end{equation}
\end{remark}

Next, we analyze how the $l^p(\Z)$ norm of $u$ changes under the Fourier rearrangement. For that, we will have to introduce the notion of symmetric rearrangement of a subset of $\Z$. We start by defining a labeling of $\Z$: Let $\epsilon: \Z \rightarrow \N$ be defined as 
\begin{align*}
    \epsilon(n) := 
    \begin{cases}
    &2n \hspace{37pt} \text{if} \hspace{9pt} n > 0\\
    &1-2n \hspace{19pt} \text{if} \hspace{9pt} n \leq 0
    \end{cases}
\end{align*}
Let $E$ be a finite subset of $\Z$ of size $k$. Then the \emph{symmetric rearrangement} $E^*$ of the set $E$ is defined as the first $k$ elements of $\Z$ with respect to labeling $\epsilon$.

\begin{lemma}\label{lem3.7}
Let $E$ be a finite subset of $\Z$ and $E^*$ be its symmetric rearrangement, then there exists a constant $c$(independent of $u$ and $E$) such that 
\begin{equation}\label{3.17}
    \sum_{n \in E} |u|^2 \leq c \sum_{n \in E^*} |u^\#|^2,
\end{equation}
for all $u \in l^2(\Z)$.
\end{lemma}

\begin{remark}
The proof of lemma \ref{lem3.7} is an adaptation of the proof of theorem 1 in \cite{montgomery1976note}. Similar ideas were used in \cite{rupert} to prove inequalities of the type \eqref{3.17} in the context of Fourier rearrangement on $\R^d$.  
\end{remark}

\begin{proof}
Let $\tau$ be a fixed non-negative number, which we will choose later. Define
\begin{equation}
    F(u_1) := \chi_{\{|F(u)| \geq \tau\}} F(u) \hspace{5pt}  \text{and} \hspace{5pt} F(u_2) := \chi_{\{|F(u)| < \tau\}} F(u).
\end{equation}
Clearly we have $u = u_1 + u_2$ and 
\begin{align*}
    \sum_{n \in E}|u|^2 \leq 2\sum_{n \in E}|u_1|^2 + 2 \sum_{n \in E}|u_2|^2.
\end{align*}
Consider
\begin{equation}\label{3.19}
    \sum_{n \in E}|u_1|^2 \leq ||u_1||_\infty^2 |E| \leq |E|(2\pi)^{-\frac{1}{2}} \Big(\int_{\{|F(u)| \geq \tau \}} |F(u)| dx\Big)^2, 
\end{equation}
and 
\begin{equation}\label{3.20}
    \sum_{n \in E}|u_2|^2 \leq \sum_{n \in \Z} |u_2|^2 = \int_{\{|F(u)| < \tau\}} |F(u)|^2 dx.
\end{equation}
It is easy to see that the terms on the RHS of \eqref{3.19} and \eqref{3.20} are invariant under symmetric decreasing rearrangement, i.e., 
\begin{align*}
    \int_{\{|F(u)| \geq \tau \}} |F(u)| dx = \int_{\{|F(u^\#)| \geq \tau \}} |F(u^\#)| dx  \hspace{5pt} \text{and} \hspace{5pt} \int_{\{|F(u)| < \tau\}} |F(u)|^2 dx =\int_{\{|F(u^\#)| < \tau\}} |F(u^\#)|^2 dx.
\end{align*}

Let $v := u^\#$, then it only remains to prove the following inequalities:
\begin{equation}\label{3.21}
    \sum_{n \in E^*}|v|^2 \gtrsim \int_{\{|F(v)| < \tau\}} |F(v)|^2 dx,
\end{equation}
and
\begin{equation}\label{3.22}
    \sum_{n \in E^*}|v|^2 \gtrsim |E| \Big(\int_{\{|F(v)| \geq \tau \}} |F(v)| dx\Big)^2. 
\end{equation}
First we prove inequality \eqref{3.21}. Let $R := (|E|-1)/2$ if $|E|$ is an odd number and $R := (|E|-2)/2$ if $|E|$ is an even number, and let $K(n)$ := max$(0, 1- |n|/R)$. Then 
\begin{align*}
    F(k)(x) = \frac{1}{R} \frac{\sin^2(Rx/2)}{\sin^2(x/2)}.
\end{align*}
It is easy to see that $F(k)(x) \geq 4R/\pi^2$ for $ |x| \leq \pi/R$. Now choose $\tau$ = inf $\{F(v)(x):  |x| \leq \pi/R\}$ and consider,
\begin{align*}
    \sum_{n \in E^*} |v|^2 \geq \sum_{n \in \Z} K(n)|v|^2 &= \int_\R \int_\R F(v)(x) F(k)(x-y) F(v)(y) dy dx \\
    & \geq \int_{x > \pi/R} \int_{x-\pi/R \leq y \leq x} F(v)(x) F(k)(x-y) F(v)(y) dy dx\\
    & \geq 2/\pi \int_{|x| > \pi/R} |F(v)(x)|^2 dx \geq 2/\pi \int_{\{|F(v)| < \tau\}} |F(v)(x)|^2 dx. 
\end{align*}
Next we prove inequality \eqref{3.22}. Consider,
\begin{align*}
   \sum_{n \in E^*} |v|^2 &\geq \int_\R \int_\R F(v)(x) F(k)(x-y) F(v)(y) dy dx \\
   & \geq 4R/\pi^2 \int_{0 \leq x \leq \pi/R } \int_{0 \leq y \leq \pi/R} F(v)(x)F(v)(y) dy dx \\
   & \geq R/\pi^2 \Big( \int_{|x| \leq \pi/R} F(v) dx \Big)^2  \geq \frac{1}{4\pi^2} |E| \Big( \int_{\{|F(u)| \geq \tau\}} F(v) dx \Big)^2.
\end{align*}
This proves inequalities \eqref{3.21} and \eqref{3.22}, thereby completing the proof of \eqref{3.17}.

\end{proof}

\begin{theorem}
Let $\varphi: [0, \infty) \rightarrow \R$ be a non-negative convex function with $\varphi(0) = \varphi'(0)= 0$. Then, for all $u \in l^2(\Z)$ we have
\begin{equation}\label{3.23}
    \sum_{n \in \Z}\varphi(|u|^2) \leq  \sum_{n \in \Z} \varphi(c|u^\#|^2), 
\end{equation}
where $c$ is the constant in \eqref{3.17}.
\end{theorem}

\begin{proof}
Consider,
\begin{align*}
    \phi(s) = \int_{0}^s \phi'(t) dt &= - \int_0^s \phi'(t)(s-t)'dt\\
    &= \int_0^s \phi''(t)(s-t)dt = \int_0^\infty     \phi''(t)(s-t)_{+} dt.
\end{align*}
Choosing $s= |u|^2$ and summing both sides we obtain,
\begin{equation}\label{3.24}
    \sum_{n \in \mathbb{Z}}\phi(|u|^2) = \int_{\mathbb{R}} \sum_{n \in \mathbb{Z}}(|u|^2-t)_{+} \phi''(t) dt.
\end{equation}
Applying \eqref{3.17} with $E = \{|u|^2 >t\}$ we get, 
\begin{align*}
    \sum_{n \in \Z}(|u|^2-t)_+ = \sum_{n \in \{|u|^2 >t\}} (|u|^2 -t) \leq \sum_{ n \in \{|u|^2 > t\}^*} (c|u^\#|^2-t) \leq \sum_{n \in \Z}(c|u^\#|^2-t)_+. 
\end{align*}
Plugging the above estimate in \eqref{3.24} completes the proof.

\end{proof}

\begin{corollary}\label{cor3.10}
For $p > 2$, there exists a constant $c(p)$ such that 
\begin{equation}\label{3.25}
    \sum_{n \in \Z}|u|^p \leq c(p) \sum_{n \in \Z} |u^\#|^p.
\end{equation}
\end{corollary}

\begin{proof}
Apply \eqref{3.23} with $\varphi(x):= x^{p/2}$.
\end{proof}

\begin{remark}
We will end this section with an \emph{open problem}. The estimates done above are very crude, and the constant $c(p)$ in \eqref{3.25} is an exponential function of $p$. It would be a worthwhile effort to find the sharp constant in \eqref{3.25}.
\end{remark}

\section{Symmetric-Decreasing rearrangement}\label{sec:symmetricdecreasing}
In this section, we will combine two rearrangements defined in the above two sections to construct a rearrangement which will be both radial and decreasing. More precisely, \\

Let $u \in l^2(\mathbb{Z})$. We define 
\begin{equation}
    v(n) := \widetilde{\Big(u^\#|_{\Z^+}\Big)}. 
\end{equation}
Then the $\emph{symmetric-decreasing rearrangement}$ $u^*$ of the function $u$ is defined as 
\begin{equation}
    u^*(n) := 
    \begin{cases}
    v(n) \hspace{28pt} \text{if} \hspace{5pt} n \geq 0\\
    v(-n) \hspace{19pt} \text{if} \hspace{5pt} n < 0\\
    \end{cases}
\end{equation}
Properties of function $u^*$:
\begin{enumerate}
    \item $u^*$ is always non-negative.\\
    \item It is easy to see that $u^*$ is radially symmetric and decreasing, that is, 
    \begin{align*}
        u^*(x) = u^*(y) \hspace{29pt} \text{if} \hspace{5pt} |x|=|y|,
    \end{align*}
    and 
    \begin{align*}
        u^*(x) \geq u^*(y) \hspace{29pt} \text{if} \hspace{5pt} |x| \leq |y|. \\
    \end{align*}
    \item $l^2(\Z)$ norm is preserved, that is,
    \begin{equation}\label{4.3}
        \sum_{n \in \Z}|u|^2 = \sum_{n \in \Z}|u^*|^2.
    \end{equation}
    This is a consequence of the fact that both the decreasing and Fourier rearrangements preserve the $l^2$ norm plus the fact that $|u^\#(n)| \leq |u^\#(0)|$.\\
    \item (\emph{Hardy-Littlewood inequality}) Let $u,v \in l^2(\Z)$, then we have
    \begin{equation}\label{4.4}
        \Big|\sum_{n \in \Z} u(n) \overline{v(n)}\Big| \leq \sum_{n \in \Z} u^*(n) v^*(n).
    \end{equation}
    This is again a consequence of Hardy-Littlewood inequality for decreasing and Fourier rearrangements along with $|u^\#(n)| \leq |u^\#(0)|$.\\
    \item $l^2(\Z)$ distance decreases under symmetric-decreasing rearrangement, 
    \begin{equation}
        \sum_{n \in \Z} |u^* - v^*|^2 \leq \sum_{n \in \Z}|u-v|^2,
    \end{equation}
    for $u, v \in l^2(\Z)$, this follows from inequality \eqref{4.3} and \eqref{4.4}.
    
\end{enumerate}
\begin{theorem}
For $p >2$, there exists a constant $c(p)$ such that,
\begin{equation}\label{4.6}
    \sum_{n \in \Z}|u|^p \leq c(p) \sum_{n \in \Z}|u^*|^p,
\end{equation}
and 
\begin{equation}\label{4.7}
     \sum_{n \in \Z}|u(n)-u(n-1)|^2 \geq \sum_{n \in \Z}|u^*(n)-u^*(n-1)|^2.
\end{equation}
for all $u \in l^2(\Z)$. 
\end{theorem}

\begin{proof}
The proof of \eqref{4.6} follows from corollary \ref{cor3.10} and the fact that $l^p(\Z)$ norm is preserved under decreasing rearrangement.  \\

Next we prove inequality \eqref{4.7}. Polya-Szeg\"{o} inequality for decreasing rearrangement \eqref{2.8} along with the radiality of $u^\#$ gives us
\begin{align*}
    \sum_{n \in \Z}|u^\#(n)-u^\#(n-1)|^2 &= 2\sum_{n \in \Z^+}|u^\#(n)-u^\#(n+1)|^2 \\
    & \geq 2\sum_{n \in \Z^+}|u^*(n)-u^*(n+1)|^2\\
    &= \sum_{n \in \Z}|u^*(n)-u^*(n-1)|^2.
\end{align*}
The above inequality along with Polya-Szeg\"{o} inequality \eqref{3.7} for Fourier rearrangement completes the proof.
\end{proof}

\begin{remark}
Let $u$ be a function which produces equality in \eqref{4.7}. Then $u$ will produce equality in the Polya-Szeg\"{o} inequality for Fourier rearrangement \eqref{3.7} and $u^\#$ will produce equality in the Polya-Szeg\"{o} inequality for decreasing rearrangement \eqref{2.8}. This would imply that $|F(u)| = F(u)^*$ and $u^\# = u^*$. We couldn't derive any useful information from these two conditions. It would be nice to understand the equality case of \eqref{4.7} in a more direct way.   

\end{remark}

\section{Applications}\label{sec:applications}
In this section, we will apply weighted Polya Szeg\"{o} inequality \eqref{2.9} to prove Hardy's inequality with the weights $n^\alpha$. In what follows, $C_c(\Z^+)$ denotes the space of finitely supported functions on $\Z^+$.

\begin{theorem}
Let $1 <\alpha \leq 2$ and $u \in C_c(\Z^+)$ with the condition that $|u(0)|^2 = max\hspace{3pt}|u|$.Then the following Hardy's inequality is true
\begin{equation}\label{5.1}
    \sum_{n=1}^\infty |u(n)-u(n-1)|^2 n^\alpha \geq (\alpha-1)^2/4 \sum_{n=1}^\infty |u|^2 n^{\alpha-2}.  
\end{equation}
\end{theorem}

\begin{proof}
Let $w(n):= n^{\alpha-2}$ on $\N$ and $w(0):=1$. Then using the Hardy-Littlewood inequality \eqref{2.5} and the weighted Polya-Szeg\"{o} inequality \eqref{2.9} for decreasing rearrangement $\widetilde{u}$ we get
\begin{align*}
    \sum_{n=1}^\infty|u|^2 n^{\alpha-2} + |u(0)|^2 =\sum_{n=0}^\infty|u|^2 w &\leq \sum_{n=0}^\infty |\widetilde{u}|^2 w(n) = \sum_{n=1}^\infty|\widetilde{u}|^2 n^{\alpha-2} + |\widetilde{u}(0)|^2,\\
    \sum_{n=1}^\infty |u(n)-u(n-1)|^2 n^\alpha &\geq \sum_{n=1}^\infty |\widetilde{u}(n)-\widetilde{u}(n-1)|^2 n^\alpha.
\end{align*}

Since $|u|$ is maximized at origin, it is sufficient to prove \eqref{5.1} for $\widetilde{u}$. Therefore from now on, we can assume that $u$ is a decreasing function. \\

Let $Lu$ be the linear interpolation of $u$ on $\R$, i.e.,
\begin{align*}
 Lu(x) := u(n) + (x-n)(u(n)-u(n-1)) \hspace{5pt} \text{for} \hspace{5pt} x \in [n-1, n]   
\end{align*}
Using the linearity of $Lu$ and the fact that $u$ is a decreasing function we can conclude that
\begin{equation}
    \sum_{n=1}^\infty|u|^2 n^{\alpha-2} \leq \int_0^\infty |Lu|^2 x^{\alpha-2} dx
\end{equation}
and
\begin{equation}
    \sum_{n=1}^\infty |u(n)-u(n-1)|^2 n^\alpha \geq \int_0^\infty |(Lu)'|^2 x^\alpha dx.
\end{equation}
We complete the proof by applying the weighted Hardy's inequality on $Lu$ on the interval $(0, \infty)$.

\end{proof}

\textbf{Acknowledgments.} We would like to thank Professor Ari Laptev for his time and various valuable discussions. We would also like to thank Rupert Frank for sharing the proof of the fact that the $L^p(\R^d)$ norm of a function increases under the Fourier rearrangement. Finally, the author wishes to thank Ashvni Narayanan for her careful proof-reading. The author is supported by the President's Ph.D. Scholarship, Imperial College London.
%%%%%%%%%%%%%%%%%%%%%%%%%%%%%%%%%%%%%%%%%%%
%%%%%%%%%%%%%%%%%%%%%%%%%%%%%%%%%%%%%%%%%%%

%\bibliographystyle{amsalpha}

\end{document}